\documentclass{amsart}

\usepackage{graphicx}
\usepackage{amssymb}
\usepackage{amsthm}
\usepackage{amsmath}
\usepackage{pstricks,pst-node}
\usepackage{longtable}
\usepackage{array}
\usepackage{ragged2e}
\swapnumbers
\newtheorem{theorem}[subsection]{Theorem}
\newtheorem{lemma}[subsection]{Lemma}

\newtheorem{cnj}[subsection]{Conjecture}

\theoremstyle{definition}
\newtheorem{definition}[subsection]{Definition}

\newcommand{\Z}{\mathbb{Z}}
\newcommand{\Q}{\mathbb{Q}}
\newcommand{\F}{\mathbb{F}}

\newcommand{\iso}{\cong}

\def\HH{{\mathcal H}}

             \def \cH {{\mathcal H}}

\def \cV {{\mathcal V}}
             \def \cW {{\mathcal W}}

\def\F{\mathbb{F}}

\def\R{\mathbb{R}}
\def\C{\mathbb{C}}
\def\Q{\mathbb{Q}}
\def\Z{\mathbb{Z}}

\DeclareMathOperator{\Frob}{Frob}
\DeclareMathOperator{\Gal}{Gal}

\newcommand{\GL}{\mathrm{GL}}
\newcommand{\SL}{\mathrm{SL}}

\usepackage{bm}
\newcommand{\ChiS}{\bm{\rho}_{2,4}} 
\newcommand{\ChiB}{\bm{\rho}_{2,3}} 

\begin{document}

\title[Mod~$2$ homology for $GL(4)$ and Galois representations]{Mod~$2$ homology for $GL(4)$ and Galois representations}

\author{Avner Ash} \address{Boston College\\ Chestnut Hill, MA 02445}
\email{Avner.Ash@bc.edu} \author{Paul E. Gunnells}
\address{University of Massachusetts Amherst\\ Amherst, MA 01003}
\email{gunnells@math.umass.edu} \author{Mark McConnell}
\address{Princeton University\\ Princeton, New Jersey 08540}
\email{markwm@princeton.edu}

\thanks{
AA wishes to thank the National Science Foundation for support of this
research through NSF grant DMS-0455240, and also the NSA through grant
H98230-09-1-0050.  This manuscript is submitted for publication with
the understanding that the United States government is authorized to
produce and distribute reprints.  PG wishes to thank the National
Science Foundation for support of this research through NSF grant
DMS-0801214 and DMS-1101640.}


\dedicatory{Dedicated to the memory of Steve Rallis}
\begin{abstract}
We extend the computations in \cite{AGM4} to find the mod~$2$ homology
in degree $1$ of a congruence subgroup $\Gamma$ of $\SL(4,\Z)$ with
coefficients in the sharbly complex, along with the action of the
Hecke algebra.  This homology group is related to the cohomology of
$\Gamma$ with $\F_2$ coefficients in the top cuspidal degree. These
computations require a modification of the algorithm to compute the
action of the Hecke operators, whose previous versions required
division by $2$.  We verify experimentally that every mod~$2$ Hecke
eigenclass found appears to have an attached Galois representation,
giving evidence for a conjecture in \cite{AGM4}.  Our method of
computation was justified in \cite{AGM5}.
\end{abstract}

\maketitle

\section{Introduction}\label{intro}

\subsection{}
This is a continuation of a series of papers \cite{AGM1, AGM2, AGM3,
  AGM4, AGM5} devoted to the computation of the cohomology of
congruence subgroups $\Gamma\subset \SL(4,\Z)$ with constant
coefficients, together with the action of the Hecke operators on the
cohomology.  We also investigated the representations of the absolute
Galois group of $\Q$ that appear to be attached to Hecke eigenclasses
in the cohomology.  The papers \cite{AGM1, AGM2, AGM3} deal with
complex coefficients, while \cite{AGM4, AGM5} deal with coefficients
in a prime finite field $\F_p$, with $p$ odd.  The current paper takes
$p=2$.  We concentrate on $H^5(\Gamma)$ because on the one hand, $H^5$
supports cuspidal cohomology with $\C$-coefficients, and on the other
hand it is only one degree below the virtual cohomological dimension
of $\Gamma$ and therefore amenable to an algorithm due to one of us
(PG) for computing Hecke operators \cite{experimental}.  Our next
project will be to rewrite our code to deal with finite-dimensional
twisted coefficients, which should lead to more interesting examples
of attached Galois representations, aimed at testing the
generalization of Serre's conjecture found in \cite{ADP} and in
\cite{Herzig}.

As explained in \cite{AGM4}, when $p>5$, the $\C$- and mod~$p$-betti
numbers coincide.  In this case we can compute the cohomology in terms
of the Steinberg module and the sharbly complex, which is what we did.
Namely, $H^5(\Gamma,K)\approx H_1(\Gamma,St\otimes K)\approx
H_1(\Gamma,Sh_{\bullet }\otimes K)$.  Here, $K=\C$ or $\F_p$, $\Gamma \subset
\SL(4,\Z)$ is a congruence subgroup, $St$ denotes the Steinberg
module, and $Sh_{\bullet }$ the sharbly complex, whose definitions are recalled
in Section~\ref{sh} below.

What we actually compute is the homology valued in the sharbly
complex.  In theory we could compute the mod~$p$ cohomology using a
spectral sequence similar to that used by Soul\'e for $\SL(3,\Z)$ in
\cite{Soule}.  However, even if we carried out this arduous task, we
do not know how to compute the Hecke action on the resulting
cohomology.

The method we use to compute $H_1(\Gamma,Sh_{\bullet }\otimes \F_2)$ is the same
as in the previous papers.  However, the algorithm in
\cite{experimental} for the computation of the Hecke action had
required division by~$2$, which prevented our treatment of mod~$2$
coefficients.  Following a suggestion of Dan Yasaki, we overcome that
problem in this paper.

\subsection{}
The mod~$2$ homology is especially interesting for two reasons.  One
is that there are many more mod~$2$ classes than exist for odd primes,
so there is more opportunity for testing conjectures and studying
phenomenology.  The other is that every mod~$2$ Galois representation
is odd and therefore again there are more possibilities for
investigating the Serre-type conjectures.

By Theorem 13 of \cite{AGM5}, the Hecke eigenvalue data we compile
gives us parts of Hecke eigenpackets ocurring in the sharbly homology
of $\Gamma_0(N)$, for various~$N$.  Therefore we can test
Conjecture~5(d) of \cite{AGM4}, which asserts the existence of a
Galois representation unramified outside $2N$ associated to each such
eigenpacket.  We do this by searching for the Galois representation
using a computer program described in Section~\ref{gal}.  There can
easily be more than one Galois representation that fits our data for
any given Hecke eigenclass, because we have only computed a few Hecke
operators at each level (because of time and space constraints).  Our
Galois finder searches for the ``simplest'' Galois representation that
fits our data in each case.  We use the supply of characters and
2-dimensional representations coming from classical modular forms of
weights~2, 3 and~4.  In no case do we fail to find a match, using just
reducible representations made out of these blocks.  An explanation of
why we use just these weights appears in Section~\ref{gal}.

Although we stop searching when we have found one Galois
representation that appears to be attached to a given Hecke
eigenpacket, we know by the Brauer-Nesbitt Theorem that up to
semisimplification there can be at most one Galois representation that
is truly attached.  This Galois representation might be describable in
many different ways using characters and classical cuspforms, because
such things can be congruent modulo a prime above~2.  Of course, we
would expect more complicated and even irreducible 4-dimensional
representations to be needed if we could compute for much larger
levels and more Hecke operators.  But at least in this small way we
find evidence both for Conjecture~5(d) of \cite{AGM4} and of the
correctness of our computations.

\subsection{}\label{ss:1.3}

As we have said, we compute Hecke eigenvalues in the sharbly homology,
not in the group homology.   When the coefficients are $\F_2$ (as in
this paper), the relationship between these two homology theories is
rather obscure, as explained in \cite{AGM4}. We might say that the
sharbly homology is more ``geometric" than the group homology, closely
related as the former is to the Voronoi decomposition of the cone of
$n$-dimensional quadratic forms. We believe that the sharbly homology
is an interesting Hecke module in its own right, conjecturally
possessing attached Galois representations.

We have no idea at present how one might prove in general that Galois
representations are attached to \ Hecke eigenclasses in the sharbly
homology.  However, we believe that any naturally occurring Hecke
module arising from the "geometry" of $\GL(n,\Z)$ ought to have
attached Galois representations, and that is what we are testing in
this paper in the case of the mod 2 sharbly homology.

When looking experimentally for apparently attached Galois
representations, naturally we can only look for Galois representations
with small image, ones we can get our hands on. Of course we would be
happy to find Hecke eigenclasses where the conjecturally attached
Galois representations have large image, especially irreducible or
non-selfdual image.  But if, for example, the predicted image were
$\GL(4,\F_2)$, it would be unlikely that we could find this Galois
representation by searching for an appropriate polynomial.

On the other hand, experience suggests that when the level is small,
the attached Galois representation will be reducible, and the smaller
the level, the more reducible.  This is borne out by the data in this
paper, and it allows our ``Galois finder" to succeed.  Of course, if
we could compute much larger levels we would expect more interesting
Galois representations to be attached, although we would probably be
unable to find them, unless they were lifts from a smaller group,
e.g.~an orthogonal or a symplectic group. In this sense, failure of
our Galois finder is to be expected in the case of sufficiently large
level, but for the small levels studied in this paper, we expected and
indeed found that our Galois finder, which only deploys 1- and
2-dimensional Galois representations, always succeeded.  The
limitations on the size of the level $N$ and prime $\ell$ in the Hecke
operators $T(\ell,k)$ that we compute come purely from limitations on
computer speed and memory size. The speed limits the Hecke operators,
because the number of single cosets in $T(\ell,k)$ grows like $\ell^3$
or $\ell^4$, depending on $k$.  The size of the memory limits $N$
because the number of rows and columns of the matrices on which we
have to perform row reduction to compute the sharbly homology grows
like $N^3$. For the largest $N$ and $\ell$ occurring in this paper,
our computations were using up all available memory and could take a
full day of CPU time to compute a single $T(\ell,k)$.

Our Galois finder has two phases. In the first phase, we find attached
Galois representations which are sums of 4 characters or sums of 2
characters and a Galois representation attached to a classical cusp form
of weight 2 or 4 (modulo 2). These weights are natural in view of the
geometry of the Borel-Serre boundary. (Note however that no one has
yet carried out a complete explicit computation of the homology of the
Borel-Serre boundary for general congruence subgroups $\GL(n,\Z)$ for any
$n>3$.)

The other phase of the Galois finder is similar except that it uses
cuspforms of weights 2 and 3. These weights can be predicted by means
of Serre's conjecture for $\GL(2)$, and thus rely on the structure of
the set of classical modular forms mod 2. See Section 4 for more
discussion of these matters. It would be very interesting to have a
more geometric or automorphic interpretation of the weight 3 forms
that appear, perhaps some kind of mod 2 endoscopic or other functorial
lifting.  But we do not know of any.

\subsection{} After the completion of this paper, we learned of the
recent remarkable work of Peter Scholze \cite{scholze}, which among
other things attaches Galois representations to Hecke eigenclasses in
the mod $p$ cohomology of locally symmetric spaces.  This result
builds on results of Harris--Lan--Taylor--Thorne \cite{hltt}, which
attaches Galois representations to Hecke eigenclasses in the
characteristic $0$ cohomology of locally symmetric spaces.  It is
quite likely that Scholze's results will imply that Galois
representations are attached to Hecke eigenclasses in the mod $p$
cohomology of congruence subgroups $\Gamma $ of $\GL(n,\Z)$, for all $p$ and
$\Gamma $, although at the moment the necessary results seem to be
conditional on stabilization of the twisted trace formula.

\subsection{} We now give a guide to the paper.  In Section~\ref{sh}
we recall the definitions of the Steinberg module, the sharbly
complex, and the concept of attached Galois representation.  We state
the conjecture of \cite{AGM4} that asserts the existence of attached
Galois representations to Hecke eigenclasses in the sharbly homology.

In Section~\ref{gu} we describe what we actually compute, namely
certain Hecke eigenclasses in the sharbly homology in degree $1$. We
use the Voronoi complex.  We describe how the sharbly homology is
calculated as a Hecke module, with reference to our earlier papers for
details.  Then we explain what modifications we made to the Hecke 
algorithm to allow us to work with $\F_2$-coefficients.

In Section~\ref{gal} we describe our Galois representation finder.
Because there are so many mod~$2$ homology classes, we had to automate
the process of finding candidates for the conjecturally attached
Galois representations.

In Section~\ref{res} we give our results.  We give the level $N$ of
$\Gamma$, the dimension of $H_1(\Gamma,Sh_{\bullet }\otimes \F_2)$,
and an enumeration of packets of Hecke eigenvalues.  For each packet,
we give the dimension of its simultaneous eigenspace and a Galois
representation that appears to be attached to the packet.

We thank Dan Yasaki for conversations that greatly helped this project
at the start.  We thank Kevin Buzzard for very helpful correspondence,
particularly in regard to~(\ref{buzztalk}).  We thank the referee for
raising the issues in subsection \ref{ss:1.3}.

\section{The Steinberg module and the sharbly complex}\label{sh}

\subsection{}
Let $n\geqslant 2$ and let $\Q^n$ denote the vector space of
$n$-dimensional row vectors.

\begin{definition}
The \emph{Sharbly complex} $Sh_{\bullet } $ is the complex of $\Z
\GL(n,\Q)$-modules defined as follows.  As an abelian group, $Sh_{k}$
is generated by symbols $[v_1,\dots,v_{n+k}]$, where the $v_i$ are
nonzero vectors in $\Q^n$, modulo the submodule generated by the
following relations:

(i) $[v_{\sigma
(1)},\dots,v_{\sigma(n+k)}]-(-1)^\sigma[v_1,\dots,v_{n+k}]$ for all
permutations $\sigma$;

(ii) $[v_1,\dots,v_{n+k}]$ if $v_1,\dots,v_{n+k}$ do not span all
of $\Q^n$; and

(iii) $[v_1,\dots,v_{n+k}]-[av_1,v_{2},\dots,v_{n+k}]$ for all $a\in
\Q^\times$.

\noindent The boundary map $\partial \colon Sh_{k} \rightarrow
Sh_{k-1} $ is given by
\[ 
\partial([v_1,\dots,v_{n+k}])=
\sum_{i=1}^{n+k} (-1)^i [v_1,\dots,\widehat{v_i},\dots v_{n+k}],
\]
where as usual $\widehat{v_i}$ means to delete $v_{i}$.
\end{definition}

The sharbly complex 
$$
\dots\to Sh_i\to Sh_{i-1}\to \dots \to Sh_1\to Sh_0
$$
is an exact sequence of $\GL(n,\Q)$-modules.  We may define the
Steinberg module $St$ as the cokernel of $\partial\colon Sh_1\to Sh_0$
(cf.~\cite[Theorem 5]{AGM5}).

Of course, all these objects depend on $n$, which we suppress from the
notation, since we will later only work with $n=4$.

Let $\Gamma$ be a congruence subgroup of $\SL(n,\Z)$.

\begin{definition} Let $M$ be a right $\Gamma$-module, concentrated in degree $0$.
The \emph{sharbly homology} of $\Gamma$ with coefficients in
$M$ is defined to be $H_*(\Gamma,Sh_{\bullet} \otimes_\Z M)$, where
$\Gamma$ acts diagonally on the tensor product.
\end{definition}

If $(\Gamma,S)$ is a Hecke pair in $\GL(n,\Z)$ and $M$ is a right
$S$-module, the Hecke algebra $\cH(\Gamma,S)$ acts on the sharbly
homology since $S$ acts (diagonally) on $Sh_{\bullet} \otimes_\Z M$.

Here is a restatement of Corollary 8 of \cite{AGM5}, which shows the close
connection between the sharbly homology and the group cohomology of
$\Gamma$:

\begin{theorem}\label{S-Steinberg}
For any $\Gamma\subset \GL(n,\Z)$ and any coefficient module $M$ in
which $2$ is invertible, there is a natural isomorphism of Hecke
modules
$$
H_*(\Gamma,Sh_{\bullet} \otimes_\Z M)\to H_*(\Gamma,St \otimes_\Z M).
$$
\end{theorem}
By Borel-Serre duality \cite{B-S}, if~$\Gamma$ is torsionfree, there
is a natural isomorphism of Hecke modules
$$
H_i(\Gamma,St \otimes_\Z M)\to H^{\binom{n}{2}-i}(\Gamma, M)
$$ for all $i$.  This result can be extended to any~$\Gamma$ as long
as its torsion primes are invertible on~$M$.

In general, the sharbly homology is more mysterious.  Nevertheless, we
still expect it to have number theoretic significance, as described in
Conjecture~\ref{conj1} as follows.

\subsection{}
Let $\Gamma_0(N)$ be the subgroup of matrices in $\SL(4,\Z)$ whose
first row is congruent to $(*,0,0,0)$ modulo $N$.  Define $S_N$ to be
the subsemigroup of integral matrices in $\GL(4,\Q)$ satisfying the
same congruence condition and having positive determinant relatively
prime to $2N$.

Let $\cH(N)$ denote the $\Z$-algebra of double cosets
$\Gamma_0(N)S_N\Gamma_0(N)$.  Then $\cH(N)$ is a commutative algebra
that acts on the cohomology and homology of $\Gamma_0(N)$ with
coefficients in any $\F_2[S_N]$ module.  When a double coset is acting
on cohomology or homology, we call it a Hecke operator.  Clearly,
$\cH(N)$ contains all double cosets of the form
$\Gamma_0(N)D(\ell,k)\Gamma_0(N)$, where $\ell$ is a prime not
dividing $2N$, $0\leqslant k\leqslant m$, and
$$D(\ell,k)=\left(\begin{matrix}
1&&&&&\cr&\ddots&&&&\cr&&1&&&\cr&&&\ell&&\cr&&&&\ddots&\cr&&&&&\ell\cr\end{matrix}\right)$$
is the diagonal matrix with the first $m-k$ diagonal entries equal to
1 and the last $k$ diagonal entries equal to $\ell$.  It is known that
these double cosets generate $\cH (N)$
(cf.~\cite[Thm.~3.20]{Shimura}).  When we consider the double coset
generated by $D(\ell,k)$ as a Hecke operator, we call it $T(\ell,k)$.

We can extend $\cH(N)$ to a larger commutative algebra $\cH^*(N)$ by
adjoining the double cosets of $D(\ell,k)$ for $\ell\mid N$.  Such a
double coset, considered as a Hecke operator, is denoted $U(\ell,k)$.

Let $\overline{\F_2}$ be an algebraic closure of $\F_2$.

\begin{definition}\label{def:hp}
Let $V$ be an $\cH(N)\otimes_\Z \overline{\F_2}$-module. Suppose that
$v\in V$ is a simultaneous eigenvector for all $T(\ell,k)$ and that
$T(\ell,k)v=a(\ell,k)v$ with $a(\ell,k)\in \overline{\F_2}$ for all
prime $\ell\not | \ 2N$ and all $0\leqslant k\leqslant 4$.  If
$$\rho\colon G_\Q\to \GL(4,\overline{\F_2})$$ is a continuous representation of
$G_{\Q} = \Gal (\overline\Q/\Q)$ unramified outside $2N$, and
\begin{equation}\label{eqn:hp}
\sum_{k=0}^{4}(-1)^k\ell^{k(k-1)/2}a(\ell,k)X^k=\det(I-\rho(\Frob_\ell)X)
\end{equation}
for all $\ell\not | \ 2N$, then we say that $\rho$ is attached to $v$. 
\end{definition}  
Here, $\Frob_\ell$ refers to an arithmetic Frobenius element, so that
if $\varepsilon$ is the cyclotomic character, we have $\varepsilon
(\Frob_\ell)=\ell$.  The polynomial in~(\ref{eqn:hp}) is called the
\emph{Hecke polynomial} for~$v$ and~$\ell$.  If $\ell\mid N$, we can
still compute the left-hand side of~(\ref{eqn:hp}) and call it the
Hecke polynomial for $U(\ell,k)$, but it has no obvious bearing on the
attached Galois representation.

The following is a special case of 
\cite[Conjecture 5]{AGM4}:

\begin{cnj}\label{conj1} 
Let $N\geqslant1$.  Let $v$ be a Hecke eigenclass in
$H_{\ast}(\Gamma_0(N),Sh_{\bullet}\otimes_\Z \overline{\F_2} )$.  Then
there is attached to~$v$ a continuous representation unramified
outside $2N$,
$$\rho\colon G_\Q\to \GL(4,\overline{\F_2}).$$
\end{cnj}

\section{Computing homology and the Hecke action mod~$2$}\label{gu}

\subsection{}
As explained in Sections 5 and 6 of \cite{AGM5}, we compute the Hecke
operators acting on sharbly cycles that are supported on Voronoi
sharblies.  Theorem 13 of \cite{AGM5} guarantees that the packets of
Hecke eignvalues we compute do occur on eigenclasses in
$H_{1}(\Gamma_0(N),Sh_{\bullet}\otimes_\Z \overline{\F_2})$.  In this
section, we recall results from \cite{AGM5} and explain how they are
modified to work with $\F_2$ coefficients.

The sharbly complex is not finitely generated as a $\Z
\SL(n,\Z)$-module, which makes it difficult to use in practice to
compute homology.  To get a finite complex to compute $H_{1}$, we use
the Voronoi complex.  We refer to \cite[Section 5]{AGM5}) for any
unexplained notation in what follows.

Let $X_{n}^{0}\subset \R^{\binom{n+1}{2}}$ be the convex cone of
positive-definite real quadratic forms in $n$-variables.  This has a
partial (Satake) compactification $(X_{n}^{0})^{*}$ obtained by
adjoining rational boundary components, which is itself a convex cone.
The space $(X_{n}^{0})^{*}$ can be partitioned into cones $\sigma
=\sigma (x_{1},\dotsc ,x_{m})$, called \emph{Voronoi cones}, where the
$x_{i}$ are contained in certain subsets of nonzero vectors from
$\Z^{n}$.  (We write elements of $\Z^{n}$ as row vectors, as we did in
Section~\ref{sh} for $\Q^{n}$.)  The cones are built as follows: each
nonzero $x_{i}\in \Z^{n}$ determines a rank $1$ quadratic form $q
(x_{i}) = {}^{t}x_{i} x_{i}\in (X_{n}^{0})^{*}$.  Let $\Pi$ be the
closed convex hull of the points $\{ q (x)\mid x\in \Z^{n}, x\not =
0\}$.  Then each of the proper faces of $\Pi$ is a polytope, and the
$\sigma$s are exactly the cones on these polytopes.  The indexing sets
are constructed in the obvious way: if $\sigma$ is the cone on
$F\subset \Pi$, and $F$ has distinct vertices $q (x_{1}),\dotsc ,q
(x_{m})$, then the indexing set is $\{\pm x_{1},\dotsc ,\pm x_{m} \}$.
We let $\Sigma$ denote the set of all Voronoi cones.

Let $X^{*}_{n}$ be the quotient of $(X^{0}_{n})^{*}$ by homotheties.
The images of the Voronoi cones are cells in $X^{*}_{n}$.  Let $\Z
V_{\bullet}$ be the oriented chain complex on these cells, graded by
dimension, and let $\Z \partial V_{\bullet}$ be the subcomplex
generated by those cells that do not meet the interior of $X^{*}_{n}$
(i.e., the image in $X^{*}_{n}$ of the positive-definite cone).  The
\emph{Voronoi complex} is then defined to be $\cV_{\bullet} = \Z
V_{\bullet}/\Z \partial V_{\bullet}$.  For our purposes, it is
convenient to reindex $\cV_{\bullet}$ by introducing the complex
$\cW_{\bullet}$, where $\cW_{k} = \cV_{n+k-1}$.  The results of
\cite{AGM5} prove that if $n\leqslant 4$, both $\cW_{\bullet}$ and
$Sh_{\bullet}$ give resolutions of the Steinberg module.  In
particular, let $\Gamma =\Gamma_{0} (N)$.  If $M$ is a
$\Z[\Gamma]$-module such that the the order of all torsion elements in
$\Gamma$ is invertible, then $H_{*} (\Gamma , \cW_{\bullet }
\otimes_{\Z}M)\approx H_{*} (\Gamma , Sh_{\bullet}\otimes_{\Z}M)$, and
furthermore by Borel--Serre duality are isomorphic (after reindexing)
to $H^{*} (\Gamma , M)$.  These two complexes can be related as
follows in our case of interest: when $n=4$, every Voronoi cell of
dimension $\leqslant 5$ is a simplex.  Thus for $0\leqslant k
\leqslant 2$, we can define a map of $\Z[\SL(4,\Z)]$-modules
$$
\theta_k\colon \cW_{k} \to Sh_{k}
$$
that takes the Voronoi cell $\sigma(v_1,\dots,v_{k+4})$ to
$\theta_k((v_1,\dots,v_{k+4})):=[v_1,\dots,v_{k+4}]$.  This allows us
to realize Voronoi cycles in these degrees in the sharbly complex.

In the current setting, in which $M\iso \overline{\F_2}$ with trivial
$\Gamma$-action, all torsion orders in $\Gamma$ are of course not
invertible in $M$.  Hence what we actually compute is more subtle.
Let $K$ denote either $\cW_{\bullet } $ or $Sh_{\bullet }$.  It is
necessary to distinguish between $H_{*}(K) = H_*(\Gamma ,K\otimes_\Z M)$,
i.e.~the homology of $\Gamma$ with coefficients in the complex
$K\otimes_\Z M$ and $H_1(K\otimes_\Gamma M)$, which is the homology of
the complex $K\otimes_\Gamma M$ and which is the bottom line of a
spectral sequence that computes $H_{*}(K)$.  The Hecke algebra $\HH$ acts
on both of these homologies when $K=Sh_{\bullet }$, and the spectral
sequence just mentioned is $\HH$-equivariant.

Thus our computation begins by computing a basis $\{x_{i} \}$ of the
homology group $H_{1} (\cW_{\bullet }
\otimes_{\Gamma }\F_{2})$.  We then compute elements $y_{i} = \theta_{1,*}
(x_{i})\in H_{1} (Sh_{\bullet }\otimes_{\Gamma }\F_{2})$.  Let $T$ be a
Hecke operator.  We compute each Hecke translate $Ty_{i}$ and then
find a sharbly cycles $z_{i}$ such that $z_{i}=Ty_{i}$ in $H_{1}
(Sh_{\bullet}\otimes_{\Gamma}\overline{\F_2})$ and such that $z_{i}$ is in
the image of the map $\theta_{1,*}$.  The inverse images
$\theta_{1,*}^{-1} (z_{i})$ can be written as linear combinations of
the cycles $x_{i}$, which gives a matrix representing the action of
$T$ from which we can find eigenclasses and eigenvalues. 

Unfortunately, as indicated in \cite[Section 6]{AGM5}, we don't know
if the map $\theta_{1,*}$ is injective.  Thus this raises the question
of what these eigenvalues mean.  The answer is provided by Theorem 13
in \cite{AGM5}, which guarantees that if we find a cycle $v$
representing a nonzero class in $H_1(\cW \otimes_\Gamma \overline
\F_2)$ such that $\theta_1(v)T$ is homologous to $a\theta_1(v)$ in
$Sh_\bullet \otimes_\Gamma \overline{\F_2}$ (for a Hecke operator
$T$), then there exists an eigenclass in $H_1(Sh_{\bullet })$ with eigenvalue $a$
for $T$.  Hence eigenvalues we find in this way do occur in the
sharbly homology and conjecturally are associated with Galois
representations as in Conjecture~5 above.

\subsection{}
Next we turn to the actual computation of the Hecke operators.  Assume
for the moment that $\Gamma$ is torsionfree.  Let $\xi = \sum n (x)x$
be a $1$-sharbly cycle mod~$\Gamma$, where all multiplicities $n (x)$
are taken to be nonzero.  We also assume for the moment that $2$ is
invertible in the coefficient module.  As described in
\cite{experimental}, we can encode $\xi$ as a collection of $4$-tuples
$(x , n (x ), \{y \}, \{L (y ) \})$ of the following data:
\begin{enumerate}
\item The $1$-sharbly $x$ appears in $\xi = \sum n (x)x$ with
multiplicity $n (x)$.
\item $\{y \}$ is the set of $0$-sharblies appearing in the boundary
of $x$.\label{steptwo}
\item For each $0$-sharbly $y$ in \eqref{steptwo}, the matrix $L (y)$
is a \emph{lift} of $y$ to $M_{4} (\Z)$.  In other words, the rows of
the matrix $L (y)$ equal the entries of $y$, up to permutation and
scaling by $\{\pm 1 \}$.\label{stepthree}
\end{enumerate}
We further require that the lift matrices in \eqref{stepthree} are
chosen $\Gamma $-equivariantly: suppose that for $x ,x'$ in the
support of $\xi $ there exist $y$ (respectively $y'$) appearing in the
boundary of $x$ (resp., $x'$) with $y = y'\gamma$ for some $\gamma \in
\Gamma $.  Then we require $L (y)= L (y')\gamma$.  Thus we have
written $\xi$ as a collection of $1$-sharblies with multiplicities and
with extra data that reflects the cycle structure of $\xi$ mod~$\Gamma$.

The congruence groups $\Gamma$ we treat are not torsionfree in
general, and we must modify the above data.  When $\Gamma$ has
torsion, it can happen that a given $0$-sharbly $y$ is taken to itself
by a element of $\Gamma$ that reverses orientation.  In the language
of \cite[Section 3.8]{AGM1}, such Voronoi cells are
\emph{nonorientable}; in that paper and its sequels \cite{AGM2, AGM3,
  AGM4} these cells are discarded when one computes $H_{1} (\cW)$.
Unfortunately, these cells are not discardable when one computes Hecke
operators using the ideas in \cite{experimental}: after applying a
Hecke operator, such $0$-sharblies must themselves be ``reduced'' to
rewrite the Hecke translate in terms of cycles in the image of
$\theta_{1}$.

The point for the current discussion is that, when encoding $\xi$ as a
$4$-tuple, any nonorientable $0$-sharbly $y$ must effectively have
more than one lift matrix chosen for it.  In particular, if $y$ is
nonorientable then we can find an orientation-reversing $\gamma$ in
the stabilizer of $y$ with $y\gamma^{2}=y$, and we must replace the
tuple $\Phi=(x,n (x),\{y \},\{L (y) \})$ in our data with a
\emph{pair} of tuples $\Phi' , \Phi ''$.  These tuples are the same as
$\Phi$ except that (i) if $\Phi$ has multiplicity $n (x)$, then $\Phi
'$, $\Phi ''$ each have multiplicity $n (x)/2$, and (ii) if $\Phi '$
has a lift matrix $L (y)$ for $y$, then $\Phi ''$ has the lift matrix
$L (y)\gamma$ for $y$ in the same position.  

Hence we ``split'' the contribution $n (x)x$ to $\xi$ into a
contribution of two $1$-sharblies, each of multiplicity $n (x)/2$, so
that we can encode it as two $4$-tuples that can maintain the
$\Gamma$-equivariance of the data.  Of course there may be more
nonorientable $0$-sharblies in the boundary of $x$ than just $y$.  If
so we continue to split tuples as needed, dividing multiplicities by
$2$ along the way.  Since $x$ has at most $5$ $0$-sharblies in its
boundary, our original $x$ gives rise to at most $2^{5}$ tuples.

We now return to the case at hand, in which $2$ is not invertible in
the coefficients.  Clearly we cannot apply the above construction to
encode $\xi$ as a collection of tuples, since we cannot replace $n
(x)$ by $n (x)/2$ if a $0$-sharbly is taken to itself by its
stabilizer.  Fortunately we are saved by an observation of Dan Yasaki:
since $-1=1$ in the coefficients, there is no distinction between
orientable and nonorientable Voronoi cells!  All Voronoi cells are
orientable; none are discarded when one builds the complex
$\cW_{\bullet}$.  The consequence is that a sharbly chain \emph{never}
becomes a cycle mod~$\Gamma$ because of orientation-reversing
self-maps on $0$-sharbles in its boundary.  Hence we never have to
divide by $2$ in building the tuples $\Phi$ to encode $\xi$.

\section{Finding attached Galois representations}\label{gal}

\subsection{}
Suppose we have a finite-dimensional $\F_2$-vector space~$V$ with a
Hecke action.  We now describe how we find Galois representations that
are conjecturally attached to Hecke eigenvectors in $V \otimes_{\F_2}
\overline{\F_2}$.  Our Galois representation finder is a Python script
built on the mathematical software package Sage~\cite{sage}.

\subsection{} We start by using the algorithm in
\cite{experimental} to compute explicitly the Hecke operators $T(\ell,
k)$ for $k=1, 2, 3$ and for~$\ell$ ranging through a set~$L$ of small
odd primes.  The operator is $U(\ell, k)$ rather than $T(\ell, k)$ if
$\ell \mid N$.  The~$L$ we use depends on~$N$ as in
Table~\ref{table_L}.  We use a larger~$L$ when~$N$ is smaller, because
the computations are faster for smaller~$N$.

\begin{table}
\caption{\label{table_L} We compute $T(\ell,k)$ and $U(\ell,k)$ at level~$N$ for the~$\ell$ shown in~$L$.}
\begin{tabular}{|l|l|}
\hline
$N$ & $L$ \\
\hline
3--10, 17 & $\{3,5,7,11,13\}$ \\
11 & $\{3,5,7,11,13,17\}$ \\
13 & $\{3,5,7,11\}$ \\
other & $\{3,5,7\}$ \\
\hline
\end{tabular}
\end{table}

Let $\F$ be the field generated over $\F_2$ by the eigenvalues of the
Hecke operators we have computed.  $\F$ is a finite extension of
$\F_2$.  We replace~$V$ with its extension of scalars $V\otimes_{\F_2}
\F$ for the rest of the discussion.

For each operator we have computed, we decompose~$V$ into
eigen\-spaces under that operator.  Then we take the common refinement
of all the decompositions.  In other words, let $E$ have the form
$\bigcap_{(\ell, k)} E_{\ell, k}$, where $E_{\ell, k}$ is any one of
the eigenspaces for the operator at $(\ell, k)$, and the intersection
is over all $\ell\in L$ and $k=1,2,3$.  We find all the non-zero~$E$
of this form, and call each a \emph{simultaneous eigen\-space}.  The
$E$'s are pairwise disjoint, and together they span a subspace of~$V$.
By construction, the Hecke eigenvalues $a(\ell, k)$ are constant on
each~$E$ and characterize it.  The function $(\ell, k) \mapsto a(\ell,
k)$ is the \emph{Hecke eigenpacket} of~$E$.

To a simultaneous eigenspace~$E$ we now attach a family of
polynomials.  Let $L' = \{\ell \in L \mid\, \ell \nmid N\}$.

\begin{definition} The \emph{polynomial system} $\mathcal{F}(E)$ is the
mapping that sends $\ell \in L'$ to the Hecke polynomial with
eigenvalues $a(\ell, k)$ defined in~(\ref{eqn:hp}).
\end{definition}

The Hecke polynomials have coefficients in the field of
eigenvalues~$\F$, but they do not necessarily split into linear
factors over that field.  We enlarge~$\F$ if necessary so that all the
Hecke polynomials for $\ell\in L'$ split into linear factors in
$\F[X]$, and again we replace~$V$ with its extension of scalars
$V\otimes_{\F_2} \F$.  The largest $\F$ we have had to work with is
$\F_{64}$ at level $N=59$, a very small field from the computational
standpoint.

\subsection{} We will be using various Galois representations~$\rho$ that
have been defined classically.  Each $\rho$ is a continuous,
semisimple representation of $G_{\Q}$ unramified outside $2N$.  It
takes values in $\GL(n', \F)$ for $n' = 1$ or~2, where $\F$ is
the particular finite extension of $\F_2$ described above.  The
characteristic polynomial of Frobenius for~$\rho$ is known and is of
degree~$n'$ for each $\ell \nmid 2N$.

\begin{definition} The \emph{polynomial system} $\mathcal{F}(\rho)$ is the
mapping that sends $\ell \in L'$ to the characteristic polynomial of
Frobenius for~$\rho$ at~$\ell$.
\end{definition}

Before we say which~$\rho$ we consider, let us describe how we
conjecturally attach a sum of $\rho$'s to a simultaneous
eigenspace~$E$.  Say that $\mathcal{F}(\rho)$ \emph{divides}
$\mathcal{F}(E)$ if, for each $\ell\in L'$, the polynomial at~$\ell$
for $\rho$ divides the polynomial at~$\ell$ for~$E$.  When one
polynomial system divides another, define the \emph{quotient system}
in the obvious way.

For a given~$E$, let $\mathcal{F} = \mathcal{F}(E)$ be its polynomial
system.  We run through a list of Galois representations~$\rho$ in
some fixed order.  The first time we find a $\rho$ (call it $\rho_1$)
whose system divides $\mathcal{F}$, we replace $\mathcal{F}$ by the
quotient system.  If the system for~$\rho_1$ divides $\mathcal{F}$
more than once (say $n_1$ times), we take the quotient $n_1$ times.
After that, we continue running through the rest of the $\rho$'s in
our fixed order.  When we find a $\rho_2$ whose system divides the new
$\mathcal{F}$, say $n_2$ times, we again replace $\mathcal{F}$ with
the quotient system.  We stop with success when $\mathcal{F}$ becomes
the trivial system, meaning all polynomials have degree zero.  We stop
with failure when we run out of $\rho$'s before $\mathcal{F}$ becomes
trivial.  In the successful cases, we say that the \emph{Galois
  representation apparently attached to}~$E$ is
$$
\rho_1^{\oplus n_1} \oplus \rho_2^{\oplus n_2} \oplus \cdots .
$$
The word ``apparently'' means that this Galois representation matches
our Hecke data as far as our data extends.

\subsection{} Now we describe the Galois representations~$\rho$ we
use.  We have two different lists of Galois representations, $\ChiS$
and $\ChiB$.  With either list, we always successfully find a Galois
representation that is apparently attached to one of our simultaneous
eigenspaces~$E$.  Specific results are in
Section~\ref{res}. The lists $\ChiS$ and $\ChiB$ are ordered, and the
order matters in the following sense.  When we split the first
representation, $\rho_1$, off of~$E$, we want $\rho_1$ to be as simple
as possible.  $\rho_2$ should be the second simplest, and so on.

In this subsection, we define $\ChiS$ and $\ChiB$.
In~(\ref{ChiReasons})--(\ref{exponent2}), we give the motivation behind
the definitions.

$\ChiS$ begins with a list of one-dimensional Galois
representations~$\chi$.  These are Dirichlet characters with value
in~$\F$, which we identify with one-dimensional representations as
usual.  Let~$M$ be the odd part of~$N$.  The definition is that a
Dirichlet character~$\chi$ belongs to $\ChiS$ if and only if the
conductor~$N_1$ of~$\chi$ is a divisor of~$M$.  Following the
intuition that a Dirichlet character with smaller conductor is simpler
than one with a larger conductor, we put the Dirichlet characters into
$\ChiS$ in order of increasing~$N_1$.  For instance, $\chi=1$ comes
first.  Sage's class \texttt{DirichletGroup} enumerates the~$\chi$ for
a given~$N_1$ automatically.  The characteristic polynomial of
Frobenius at~$\ell$ for~$\chi$ is $1 + \chi(\ell) X$, for all $\ell
\nmid 2N$.

After the Dirichlet characters, we put into $\ChiS$ certain Galois
representations~$\rho$ coming from classical cusp forms for congruence
subgroups of $\SL(2,\Z)$.  We emphasize that the cusp forms are in
characteristic zero, though the~$\rho$ take values in characteristic
two.  The characteristic polynomials of Frobenius for the cusp forms
are naturally defined over number fields, so, as we describe which
cusp forms we use, we must also describe how we reduce to get Galois
representations defined over~$\F$.

Let $N_1$ be a divisor of~$M$.  Let $f$ be a newform of weight 2 or 4
for $\Gamma_0(N_1)$.  The coefficients of the $q$-expansion of~$f$
generate a number field~$K_f$, with ring of integers
$\mathcal{O}_{K_f}$.  Let $\mathfrak{p}$ be a prime of~$K_f$ over~2.
If $\F$ is of high enough degree over $\F_2$, then the finite field
$\mathcal{O}_{K_f} / \mathfrak{p}$ has an embedding
$\alpha_{\mathfrak{p}}$ into $\F$.  In every case we have computed,
$\F$ is indeed large enough so that this embedding exists.  Then the
pair $(f, \mathfrak{p})$ gives rise to a Galois representation~$\rho$
into $\GL(2, \F)$, by reduction mod~$\mathfrak{p}$ composed
with $\alpha_{\mathfrak{p}}$.  For any $\ell \nmid 2N$, the
characteristic polynomial of Frobenius is $1 -
\alpha_{\mathfrak{p}}(a_\ell) X + X^2$, where $a_\ell$ is the
$\ell$-th coefficient in the $q$-expansion of~$f$.

By definition, $\ChiS$ contains the representation~$\rho$ for $(f,
\mathfrak{p})$, for all $N_1 \mid M$ and all newforms $f$ of weight 2
or 4 for $\Gamma_0(N_1)$.  The order is as follows.  The outermost
loop is over weight~2 first, then weight~4.  For a given weight, we
let $N_1$ run through the divisors of~$M$ in increasing order.  We
find the newforms~$f$ for $\Gamma_0(N_1)$ and the given weight.
Sage's class \texttt{CuspForms}, with its method \texttt{newforms},
makes this last step automatic.  For each newform~$f$, we find the
number field $K_f$.  If there is more than one~$f$ for the given
weight and~$N_1$, we sort these $f$'s by two keys; the primary key
says the degree $[K_f : \Q]$ should be increasing, and the secondary
key says that the absolute value of the discriminant of $K_f$ should
be increasing.

We now turn to the definition of $\ChiB$.  It begins with the same
Dirichlet characters as $\ChiS$, in the same order.  Next, let $N_1
\mid M$.  Let $\psi$ be a character on $\Z/N_1\Z$.  Let $f$ be a
newform of weight 2 or 3 with level~$N_1$ and nebentype
character~$\psi$.  Let $K_f$ and~$\mathfrak{p}$ be as before.  The
pair $(f, \mathfrak{p})$ gives rise to a Galois representation~$\rho$
into $\GL(2, \F)$ as above.

By definition, $\ChiB$ contains the representation~$\rho$ for $(f,
\mathfrak{p})$, for all $N_1 \mid M$, all~$\psi$, and all newforms $f$
of weight 2 or 3 of level~$N_1$ and nebentype character~$\psi$.  The
order is as follows.  The outermost loop is over weight~2 first, then
weight~3.  For a given weight, we let $N_1$ run through the divisors
of~$M$ in increasing order.  For a given~$N_1$, we run through
the~$\psi$ in the order Sage uses, which is to fix generators of the
character group and raise them to powers in lexicographic order,
starting with 0-th powers.  In particular, the trivial~$\psi$ comes
first.  We find the newforms for the given weight and~$\psi$, again
using Sage's class \texttt{CuspForms}.  For each newform~$f$, we find
the number field $K_f$, and sort the $f$'s by degree and discriminant
as before.

\subsection{}\label{ChiReasons}
The definitions of $\ChiS$ and $\ChiB$ present two different
perspectives on how Galois representations would be attached to our
homology classes.

Our construction of $\ChiS$ reflects a guess based on an analogy with
our first papers in this series, which studied homology in
characteristic zero \cite{AGM1, AGM2, AGM3}.  In them, we found that,
for small levels, all the homology appeared to be accounted for by
classes supported on the Borel-Serre boundary, and that it was always
related to Dirichlet characters and classical cuspforms of weights~2
and~4.  Although with this guess we might expect to need newforms of
even level dividing $N$, in practice we did not.

By contrast, the list $\ChiB$ reflects the conjecture found
in~\cite{ADP}.  Here we seek mod~2 Galois representations which the
conjecture would associate to a homology class of level~$N$,
weight~$k$ and trivial nebentype.  In particular, the Serre conductor
of such Galois representations would divide the odd part~$M$ of the
level~$N$.  We never get to test this for more than two-dimensional
Galois representations, not all the way to four dimensions, because we
keep splitting off the Dirichlet characters.  Our guesses for
two-dimensional representations are that they are mod~2 Galois
representations which Serre's conjecture would attach to a homology
class of level~$M$, weight~$k$ and trivial nebentype.  But we don't
have a good way to construct these mod~2 objects, except by reducing
characteristic-zero modular forms mod~2.  Kevin Buzzard tells us that
we are guaranteed to find all such two-dimensional mod~2 Galois
representations by looking at modular forms of level~$M$, weights 2 or
3, and a range of nebentypes.  In~(\ref{buzztalk}) let us explain this
guarantee.

\subsection{}\label{buzztalk} We thank Kevin Buzzard for much of the
information in this subsection.  Let $\sigma$ be a mod~2 Galois
representation which Serre's conjecture (now a theorem of
Khare-Wintenberger \cite{KW1,KW2}) would attach to a homology class of
level~$M$, weight~$k$ and trivial nebentype.  The following arguments
are valid for all~$p$ prime to~$M$, including $p=2$
\cite[Thm.~4.3]{Eidx}.  First of all, we do not have to worry about
$k=1$, because by multiplying by the Hasse invariant we can move to
$k=p=2$.  Thus we may assume $k\geqslant 2$.  Any eigenform will show
up, up to a twist, in weight at most $p+1$.  Thus for $p=2$, where
there are no twists at all, we need only compute in weights $k=2$
and~3.  All the mod~2 eigenforms lift to characteristic zero, because
$k\geqslant 2$.  Because $p\leqslant 3$, we cannot guarantee that the
nebentype character lifts to the character we expect; however, we know
that it lifts to \emph{some} character.  Since the eigenforms lift to
characteristic zero, well-known work of Deligne attaches $p$-adic
representations to them.  In turn, these reduce to mod~$p$
representations, one of which is the given~$\sigma$.

In practice, there are a large number of nebentypes~$\psi$, and often
a large number of cusp forms for a given nebentype.  We cut down on
the amount of computation as follows.  Our desired four-dimensional
Galois representation must have determinant~1.  In every case where we
use a cusp form, we have already split off two Dirichlet characters;
that is, the four-dimensional representation is $\chi_1 \oplus \chi_2
\oplus \rho_3$ where $\chi_1$, $\chi_2$ are Dirichlet characters and
$\rho_3$ is from a cusp form of nebentype~$\psi$.  Thus $\det \rho_3 =
\delta$, where we define $\delta = \det(\chi_1 \chi_2)^{-1}$.
Furthermore, $\det\rho_3 = \det\psi$.  When we construct the list
$\ChiB$ is our program, we already know~$\chi_1$ and~$\chi_2$, so we
make the list smaller by only including $\psi$ that are congruent to
$\delta$ mod~2.  

\subsection{}\label{exponent2}
Let~$\Delta$ be the group of characters $\psi : (\Z/M\Z)^\times \to
\C$ that are congruent to 1 mod~2.  The $\psi$ we need to use are the
coset of~$\Delta$ translated by~$\delta$, so we would like to
understand~$\Delta$.  Let $\mu$ be the exponent of the group
$(\Z/M\Z)^\times$.  All our~$\psi$ take values in $\Q(\zeta_\mu)$, the
cyclotomic field of $\mu$-th roots of unity, and ``mod~2'' means
modulo a prime ideal $\mathfrak{p}_\mu$ over~2 in $\Q(\zeta_\mu)$.
Let $\nu$ be the power of~2 dividing~$\mu$, and let $o$ be the odd
part, so that $\mu = \nu o$.  As usual, $\Q(\zeta_\mu)$ is the
compositum of $\Q(\zeta_\nu)$ and $\Q(\zeta_o)$, and
$\mathfrak{p}_\mu$ can be understood by studying the primes
$\mathfrak{p}_\nu$, $\mathfrak{p}_o$ over~2 in their respective
fields.

\begin{lemma}
$\Delta$ is the group of characters whose image lies in
  $\Q(\zeta_\nu)$.  Equivalently, it is the group of characters whose
  orders are pure powers of two dividing~$\nu$.
\end{lemma}

To prove the lemma, first consider~$\nu$.  In $\Q(\zeta_\nu)$, 2 is
totally ramified, and $\mathfrak{p}_\nu = (2, 1 - \zeta_\nu)$ is the
only prime over~2.  Any Dirichlet character is 1 mod~2, because
$\zeta_\nu$ and all its powers are congruent to~1
mod~$\mathfrak{p}_\nu$.  Second, consider~$o$.  For an odd prime~$q$
that divides~$o$, let $o'$ be maximal power of~$q$ that divides~$o$.
In $\Q(\zeta_{o'})$, 2 is unramified.  Under the mapping to the
residue class field, the $o'$ distinct powers of $\zeta_{o'}$ map to
$o'$ distinct values, so only the trivial power $\zeta_{o'}^0 = 1$
maps to~1 mod~2.  That is, only the trivial Dirichlet character is 1
mod $\mathfrak{p}_{o'}$.  The lemma follows from the Chinese remainder
theorem.  \qed

\section{Results}\label{res}

For the list $\ChiS$, subsection~(\ref{resultsS}) contains a table of
results for several levels~$N$.  For each level~$N$, we first give the 
overall dimension of the~$H_1$ we compute.  Each succeeding row
describes a simultaneous eigenspace~$E$.  The first two columns in the
row give the type of~$E$, a Roman numeral to be defined below,
followed by $\dim E$.

Let $\mathbf{1}$ be the trivial one-dimensional Galois representation.
Roman numeral~I means that the Galois representation apparently
attached to our Hecke eigenspace is the sum of four trivial
representations, $\mathbf{1}^{\oplus 4}$.  The symbol $\mathrm{I}_m$
means the representation is the sum of two trivial and two non-trivial
representations, $\mathbf{1} \oplus \mathbf{1} \oplus \chi_m \oplus
\bar{\chi}_m$.  The non-trivial representations go to $\F_4$ rather
than $\F_2$.  More precisely, $\chi_m$ maps $(\Z/m\Z)^\times$
surjectively to $\F_4^\times$, and $\bar{\chi}_m$ is its conjugate
under $\Gal(\F_4/\F_2)$.  These statements characterize $\chi_m$ and
$\bar{\chi}_m$ up to conjugation.

Roman numerals~II and~IV mean that the Galois representation
apparently attached to our Hecke eigenspace is the sum of two
$\mathbf{1}$'s and the Galois representation attached to a cuspidal
newform from $\ChiS$.  The newform has weight~2 or~4, respectively.
The congruence subgroup is $\Gamma_0(N)$, where~$N$ is the level where
the representation first appears in the tables.

In subsection~(\ref{resultsB}), we present data for $\ChiB$, but
list only the representations where $\ChiS$ and $\ChiB$ give different
results.  Roman numeral III stands for the sum of two $\mathbf{1}$'s
and a cuspidal newform of weight~3 from $\ChiB$.

Our tables do not give the Hecke polynomials of the $T(\ell,k)$.  This
is because the Hecke polynomials can easily be recovered from the
Galois representation.  For example, all the $T$'s for a type~I
representation have Hecke polynomial $(x+1)^4$.  The Hecke polynomials
for the~$U(\ell,k)$ are described below.  The list of~$\ell$'s used
for a given~$N$ was given in Table~\ref{table_L}.

For types~II, III, and~IV, we give details about the cusp form in the
third column of each row.  The coefficients of the $q$-expansions are
in~$\Q$ unless the number field is indicated.  We write $i =
\sqrt{-1}$ as usual.

We observe that the results for types~I, $\mathrm{I}_m$, and~II are
always the same at a given level for $\ChiS$ and $\ChiB$.  The only
differences we see are when type~IV changes to type~III.  It is
somewhat surprising that the type~II representations never change
between $\ChiS$ and $\ChiB$.  For $\ChiS$ and weight~2, we always
searched for cusp forms on $\Gamma_0(N_1)$, which means
$\Gamma_1(N_1)$ with trivial nebentype. For $\ChiB$ and weight~2, we
searched for cusp forms of all nebentypes.  The observation is that
our program produced a weight~2 cusp form for some nebentype if and
only if that nebentype was trivial.

The same cusp form can appear at the same level~$N$ for different
simultaneous eigenspaces.  This reflects the different embeddings of
the number fields into~$\F$.  For example, in the last table
of~(\ref{resultsS}), with $\ChiS$ and $N=59$, the same weight-2 cusp
form appears four times, in all four of the type~II representations.
The cusp form is defined over a quintic extension of~$\Q$.  We find
that~2 factors in the quintic field as a product $\mathfrak{p}_1
\mathfrak{p}_2$ of prime ideals, where~$\mathfrak{p}_1$ is unramified
and has residue class field $\F_8$, while~$\mathfrak{p}_2$ has
ramification index~2 and residue class field $\F_2$.  The first three
occurrences of the cusp form belong to~$\mathfrak{p}_1$.  Let
$\varpi$ be a root of $x^3 + x + 1 = 0$ in $\F_8$.  The Hecke
polynomials for the first representation are
\begin{eqnarray*}
(x+1)^2 \cdot (x^2 + \varpi^2 x + 1) &\qquad& (\ell=3) \\
(x+1)^2 \cdot (x^2 + \varpi   x + 1) &\qquad& (\ell=5) \\
(x+1)^2 \cdot (x^2 + (\varpi^2 + \varpi) x + 1) &\qquad& (\ell=7) 
\end{eqnarray*}
The Galois group $\Gal(\F_8/\F_2)$ permutes $\varpi$, $\varpi^2$, and
$\varpi^4 = \varpi+\varpi^2$ in a three-cycle.  Checking the Hecke
polynomials of the second and third Galois representations, we see
that these three representations (those with eigenspaces of
dimension~4) are permuted in a three-cycle by the Galois group.  The
fourth occurrence of the cusp form (dimension~15) is
for~$\mathfrak{p}_2$; here the coefficients of the Hecke polynomials
are down in~$\F_2$.

For a given~$N$, the sum of the dimensions of the simultaneous
eigenspaces is often less than the dimension of the full~$H_1$.  This
is because many Hecke operators, both $T(\ell,k)$ and $U(\ell,k)$,
turn out not to be semisimple.

Every level~$N$ we have computed has some representations of type~I.
A few have type $\mathrm{I}_m$.  To avoid cluttering the tables, we
list these representations here.  The notation $(N,d)$ means level~$N$
has a representation with corresponding eigenspace of dimension~$d$.
When the same~$N$ occurs in more than one pair, there are simultaneous
eigenspaces where the $T(\ell,k)$ act the same but the $U(\ell,k)$ act
differently.
\begin{itemize}
\item The type~I representations that appear to be attached to our
  data are $(3,1)$, $(4,1)$, $(5,1)$, $(6,5)$, $(7,3)$, $(8.6)$,
  $(9,1)$, $(9,4)$, $(10,7)$, $(11,1)$, $(12,19)$, $(13,1)$,
  $(14,13)$, $(15,14)$, $(16,17)$, $(17,6)$, $(18,5)$, $(18,16)$,
  $(19, 1)$, $(20,30)$, $(21,16)$, $(22, 5)$, $(23,3)$, $(24,55)$,
  $(25,1)$, $(25,9)$, $(26,7)$, $(27,1)$, $(27,3)$, $(27,4)$,
  $(28,43)$, $(29,1)$, $(30,59)$, $(31,3)$, $(32,40)$, $(33, 14)$,
  $(34,29)$, $(35,18)$, $(36, 19)$, $(36,50)$, $(37,1)$, $(38,5)$,
  $(39,21)$, $(41,8)$, $(43,10)$, $(47,3)$, $(53,1)$, $(59,1)$.
  \item We find one type $\mathrm{I}_9$ representation of dimension~2
    at level~27.
  \item We find one type $\mathrm{I}_7$ representation of dimension~4
    at level~35.
\end{itemize}

We use the operators $U(\ell,k)$ to divide up the simultaneous
eigenspaces as finely as possible, and we compute their Hecke
polynomials, but we do not consider the $U(\ell,k)$ when attaching
Galois representations.  Again, to avoid cluttering the tables with
$U(\ell,k)$ data, we summarize their Hecke polynomials here.  The
general rule is that, when $\ell$ is an odd prime dividing~$N$, the
Hecke polynomial of~$U$ is $x^4+x^3+x^2+x+1$.  We list the exceptions
in the format ($N$, $U_\ell$, $d$), which means that for all the
representations with a $d$-dimensional eigenspace we have found at
level~$N$, the operator $U(\ell,k)$ has the Hecke polynomial
described.
\begin{itemize}
\item The Hecke polynomial is $(x^2+x+1)^2$ for (9, $U_3$, 4), (18,
  $U_3$, 16), (25, $U_5$, 2 or 9), (27, $U_3$, 2 or 4) and (36, $U_3$,
  50).
\item The Hecke polynomial is $x^4+x^3+1$ for (27, $U_3$, 3).
\item Let $\omega$ be a primitive cube root of unity in $\F_4$.  At
  level $N=33$, the Hecke polynomial for $U_3$ is $x^4 + \omega x^3 +
  x^2 + \omega x + 1$ on one of the 4-dimensional eigenspaces; for
  the other 4-dimensional eigenspace, it is the polynomial's
  conjugate under $\Gal(\F_4/\F_2)$, namely $x^4 + (\omega+1) x^3 +
  x^2 + (\omega+1) x + 1$.  At level $N=39$, the same pair of
  conjugate Hecke polynomials occur for $U_3$ and the pair of
  2-dimensional eigenspace, for both types III and~IV.
\item The Hecke polynomial is $x^4+x+1$ for (33, $U_3$, 9), and also
  for (39, $U_3$, 4) for both types III and~IV.
\end{itemize}

\subsection{}\label{resultsS} Here are the results of types~II and~IV
for $\ChiS$.  (Types~I and~$\mathrm{I}_m$ were described above.)

\newcolumntype{x}[1]{>{\RaggedRight}p{#1}}

\begin{center}
\begin{longtable}{|l|l|x{5in}|}
\hline \multicolumn{3}{|l|}{\textbf{Level 11}.  Dimension 5.}
\\ \hline II & 4 & $\rho_{2,11} = q - 2q^{2} - q^{3} + 2q^{4} + q^{5}
+ O(q^{6})$ \\ 
\hline
\hline
\multicolumn{3}{|l|}{\textbf{Level 13}.  Dimension 5.} \\
\hline
IV & 2 & $\rho_{4,13} = q - 5q^{2} - 7q^{3} + 17q^{4} - 7q^{5} + O(q^{6})$ \\
\hline
\hline
\multicolumn{3}{|l|}{\textbf{Level 19}.  Dimension 9.} \\
\hline
II & 4 & $\rho_{2,19} = q - 2q^{3} - 2q^{4} + 3q^{5} + O(q^{6})$ \\
\hline
IV & 2 & $\rho_{4,19} = q - 3q^{2} - 5q^{3} + q^{4} - 12q^{5} + O(q^{6})$ \\
\hline
\pagebreak[2]\hline
\multicolumn{3}{|l|}{\textbf{Level 23}.  Dimension 12.} \\
\hline
II & 9 & $q + b_{0}q^{2} + \left(-2 b_{0} - 1\right)q^{3} + \left(-b_{0} - 1\right)q^{4} + 2 b_{0}q^{5} + O(q^{6}),$ with $b_0 = (-1 + \sqrt{5})/2$. \\
\hline
\pagebreak[2]\hline
\multicolumn{3}{|l|}{\textbf{Level 25}.  Dimension 14.} \\
\hline
IV & 2 & $q - q^{2} - 7q^{3} - 7q^{4} + O(q^{6})$ \\
\hline
\pagebreak[2]\hline
\multicolumn{3}{|l|}{\textbf{Level 26}.  Dimension 25.} \\
\hline
IV & 10 & $\rho_{4,13}$ \\
\hline
\pagebreak[2]\hline
\multicolumn{3}{|l|}{\textbf{Level 27}.  Dimension 20.} \\
\hline
II & 4 & $q - 2q^{4} + O(q^{6})$ \\
\hline
\pagebreak[2]\hline
\multicolumn{3}{|l|}{\textbf{Level 29}.  Dimension 17.} \\
\hline
II & 5 & $q + b_{0}q^{2} - b_{0}q^{3} + \left(-2 b_{0} - 1\right)q^{4} - q^{5} + O(q^{6}),$ with $b_0 = -1 + \sqrt{2}$. \\
\hline
\pagebreak[2]\hline
\multicolumn{3}{|l|}{\textbf{Level 31}.  Dimension 16.} \\
\hline
II & 9 & $q + b_{0}q^{2} - 2 b_{0}q^{3} + \left(b_{0} - 1\right)q^{4} + q^{5} + O(q^{6}),$ with $b_0 = (1 + \sqrt{5})/2$. \\
\hline
\pagebreak[2]\hline
\multicolumn{3}{|l|}{\textbf{Level 33}.  Dimension 35.} \\
\hline
II & 4 & $\rho_{2,11}$ \\
\hline
II & 4 & $\rho_{2,11}$ \\
\hline
II & 9 & $\rho_{2,11}$ \\
\hline
\pagebreak[2]\hline
\multicolumn{3}{|l|}{\textbf{Level 37}.  Dimension 21.} \\
\hline
II & 12 & $q - 2q^{2} - 3q^{3} + 2q^{4} - 2q^{5} + O(q^{6})$ \\
\hline
IV & 2 & $\rho_{4,37} = q + b_{0}q^{2} + \left(-\frac{1}{8} b_{0}^{3} - \frac{9}{8} b_{0}^{2} - \frac{13}{4} b_{0} - \frac{11}{4}\right)q^{3}+ \left(b_{0}^{2} - 8\right)q^{4}  $ \newline $+ \left(\frac{13}{8} b_{0}^{3} + \frac{85}{8} b_{0}^{2} + \frac{25}{4} b_{0} - \frac{93}{4}\right)q^{5}+ O(q^{6}),$ with $b_{0}^{4} + 6 b_{0}^{3} - b_{0}^{2} - 16 b_{0} + 6=0$. \\
\hline
IV & 2 & $\rho_{4,37}$ \\
\hline
\pagebreak[2]\hline
\multicolumn{3}{|l|}{\textbf{Level 38}.  Dimension 40.} \\
\hline
II & 15 & $\rho_{2,19}$ \\
\hline
IV & 10 & $\rho_{4,19}$ \\
\hline
\pagebreak[2]\hline
\multicolumn{3}{|l|}{\textbf{Level 39}.  Dimension 41.} \\
\hline
IV & 2 & $\rho_{4,13}$ \\
\hline
IV & 2 & $\rho_{4,13}$ \\
\hline
IV & 4 & $\rho_{4,13}$ \\
\hline
\pagebreak[2]\hline
\multicolumn{3}{|l|}{\textbf{Level 43}.  Dimension 26.} \\
\hline
IV & 2 & $\rho_{4,43} = q + b_{0}q^{2} + \left(\frac{1}{8} b_{0}^{3} + \frac{1}{8} b_{0}^{2} - \frac{7}{2} b_{0} - \frac{13}{4}\right)q^{3} + \left(b_{0}^{2} - 8\right)q^{4} + \left(-\frac{9}{8} b_{0}^{3} - \frac{33}{8} b_{0}^{2} + \frac{19}{2} b_{0} + \frac{5}{4}\right)q^{5} + O(q^{6}),$ with $b_{0}^{4} + 4 b_{0}^{3} - 9 b_{0}^{2} - 14 b_{0} + 2=0$. \\
\hline
IV & 2 & $\rho_{4,43}$ \\
\hline
\pagebreak[2]\hline
\multicolumn{3}{|l|}{\textbf{Level 47}.  Dimension 25.} \\
\hline
II & 9 & $\rho_{2,47} = q + b_{0}q^{2} + \left(b_{0}^{3} - b_{0}^{2} - 6 b_{0} + 4\right)q^{3} + \left(b_{0}^{2} - 2\right)q^{4} + \left(-4 b_{0}^{3} + 2 b_{0}^{2} + 20 b_{0} - 10\right)q^{5} + O(q^{6}),$ with $b_{0}^{4} - b_{0}^{3} - 5 b_{0}^{2} + 5 b_{0} - 1=0$. \\
\hline
II & 9 & $\rho_{2,47}$ \\
\hline
\pagebreak[2]\hline
\multicolumn{3}{|l|}{\textbf{Level 53}.  Dimension 33.} \\
\hline
II & 8 & $q - q^{2} - 3q^{3} - q^{4} + O(q^{6})$ \\
\hline
IV & 2 & $q + b_{1}q^{2} + \left(\frac{1}{14} b_{1}^{3} - \frac{3}{14} b_{1}^{2} - \frac{37}{14} b_{1} - \frac{3}{2}\right)q^{3}+\left(b_{1}^{2} - 8\right)q^{4} + \left(-\frac{5}{14} b_{1}^{3} - \frac{13}{14} b_{1}^{2} + \frac{31}{14} b_{1} + \frac{3}{2}\right)q^{5}+ O(q^{6}),$ with $b_{1}^{4} + 4 b_{1}^{3} - 16 b_{1}^{2} - 42 b_{1} + 49=0$. \\
\hline
\pagebreak[2]\hline
\multicolumn{3}{|l|}{\textbf{Level 59}.  Dimension 36.} \\
\hline
II & 4 & $\rho_{2,59} = q + b_{0}q^{2} + \left(-\frac{1}{4} b_{0}^{4} + \frac{5}{4} b_{0}^{2} - \frac{1}{2} b_{0}\right)q^{3} + \left(b_{0}^{2} - 2\right)q^{4}  $ \newline $ + \left(\frac{3}{4} b_{0}^{4} + \frac{1}{2} b_{0}^{3} - \frac{23}{4} b_{0}^{2} - 3 b_{0} + 7\right)q^{5} + O(q^{6}),$ with $b_{0}^{5} - 9 b_{0}^{3} + 2 b_{0}^{2} + 16 b_{0} - 8=0$. \\
\hline
II & 4 & $\rho_{2,59}$ \\
\hline
II & 4 & $\rho_{2,59}$ \\
\hline
II & 15 & $\rho_{2,59}$ \\
\hline
IV & 4 & $q + b_{1}q^{2} + \left(-3 b_{1} + 1\right)q^{3} + \left(b_{1} - 4\right)q^{4} + \left(3 b_{1} - 17\right)q^{5} + O(q^{6}),$ with $b_1 = (1 + \sqrt{17})/2.$ \\
\hline
\end{longtable}
\end{center}

\newpage

\subsection{}\label{resultsB} Here are the results for $\ChiB$, where they differ from $\ChiS$.

\begin{center}
\begin{longtable}{|l|l|p{5in}|}
\hline \multicolumn{3}{|l|}{\textbf{Level 13}.  Dimension 5.}
\\ \hline III & 2 & $\rho_{3,13} = q + b_{0}q^{2} +
\left(\left(i - 1\right) b_{0} - 3\right)q^{3} +
\left(\left(-2 i - 2\right) b_{0} - i\right)q^{4}+ \left(b_{0} + 3 i + 3\right)q^{5} +
O(q^{6}),$ for $\Gamma_1(13)$ with nebentype mod 13 mapping $2
\mapsto i$, with coefficients in $\mathbb{Q}(i)[b_{0}]/(b_{0}^{2} + \left(2 i + 2\right) b_{0} - 3
i)$. \\ \hline
\pagebreak[2]\hline
\multicolumn{3}{|l|}{\textbf{Level 19}.  Dimension 9.} \\
\hline
III & 2 & $\rho_{3,19} = q + b_{1}q^{2} - b_{1}q^{3} - 9q^{4} + 4q^{5}
+ O(q^{6}),$ for $\Gamma_1(19)$ with nebentype mod 19 mapping $2
\mapsto -1$, where $b_1 = \sqrt{-13}$. \\
\hline
\pagebreak[2]\hline
\multicolumn{3}{|l|}{\textbf{Level 25}.  Dimension 14.} \\
\hline
III & 2 & $q + b_{0}q^{2} + i b_{0}q^{3} - iq^{4} + O(q^{6}),$ for
$\Gamma_1(25)$ with nebentype mod 25 mapping $2 \mapsto i$, with coefficients in $\mathbb{Q}(i)[b_{0}]/(b_{0}^{2} - 3 i)$. \\
\hline
\pagebreak[2]\hline
\multicolumn{3}{|l|}{\textbf{Level 26}.  Dimension 25.} \\
\hline
III & 10 & $\rho_{3,13}$ \\
\hline
\pagebreak[2]\hline
\multicolumn{3}{|l|}{\textbf{Level 37}.  Dimension 21.} \\
\hline
III & 2 & $\rho_{3,37} = q + b_{0}q^{2} + \left(\frac{1}{4} i b_{0}^{4} + \left(\frac{1}{4} i - \frac{1}{4}\right) b_{0}^{3} + \frac{11}{4} b_{0}^{2} + \left(\frac{5}{4} i + \frac{5}{4}\right) b_{0} - 3 i\right)q^{3} + \left(b_{0}^{2} - 4 i\right)q^{4} + \left(-\frac{1}{4} i b_{0}^{5} + \left(-\frac{1}{2} i + \frac{1}{2}\right) b_{0}^{4} - \frac{13}{4} b_{0}^{3} + \left(-5 i - 5\right) b_{0}^{2} + \frac{17}{2} i b_{0} + 6 i - 6\right)q^{5} + O(q^{6}),$ \newline for $\Gamma_1(37)$ with nebentype mod 37 mapping $2 \mapsto i$, with coefficients in \newline $\mathbb{Q}(i)[b_{0}]/(b_{0}^{6} + \left(3 i + 3\right) b_{0}^{5} - 10 i b_{0}^{4} + \left(-34 i + 34\right) b_{0}^{3} - 5 b_{0}^{2} + \left(-59 i - 59\right) b_{0} - 24 i)$. \\
\hline
III & 2 & $\rho_{3,37}$ \\
\hline
\pagebreak[2]\hline
\multicolumn{3}{|l|}{\textbf{Level 38}.  Dimension 40.} \\
\hline
III & 10 & $\rho_{3,19}$ \\
\hline
\pagebreak[2]\hline
\multicolumn{3}{|l|}{\textbf{Level 39}.  Dimension 41.} \\
\hline
III & 2 & $\rho_{3,39} = q + b_{0}q^{2} + \left(\left(i - 1\right) b_{0} - 3\right)q^{3} + \left(\left(-2 i - 2\right) b_{0} - i\right)q^{4} + \left(b_{0} + 3 i + 3\right)q^{5} + O(q^{6}),$ for $\Gamma_1(13)$ with nebentype mod 13 mapping $2 \mapsto i$, with coefficients in $\mathbb{Q}(i)[b_{0}]/(b_{0}^{2} + \left(2 i + 2\right) b_{0} - 3 i)$. \\
\hline
III & 2 & $\rho_{3,39}$ \\
\hline
III & 4 & $\rho_{3,39}$ \\
\hline
\pagebreak[2]\hline
\multicolumn{3}{|l|}{\textbf{Level 43}.  Dimension 26.} \\
\hline
III & 2 & $\rho_{3,43} = q + b_{1}q^{2} + \left(-\frac{1}{4} b_{1}^{5} - \frac{15}{4} b_{1}^{3} - \frac{25}{2} b_{1}\right)q^{3} + \left(b_{1}^{2} + 4\right)q^{4} + \left(\frac{1}{4} b_{1}^{5} + \frac{11}{4} b_{1}^{3} + \frac{9}{2} b_{1}\right)q^{5} + O(q^{6}),$ for $\Gamma_1(43)$ with nebentype mod 43 mapping $3 \mapsto -1$, with coefficients in $\mathbb{Q}[b_{1}]/(b_{1}^{6} + 20 b_{1}^{4} + 121 b_{1}^{2} + 214)$. \\
\hline
III & 2 & $\rho_{3,43}$ \\
\hline
\pagebreak[2]\hline
\multicolumn{3}{|l|}{\textbf{Level 53}.  Dimension 33.} \\
\hline
III & 2 & $q + b_{0}q^{2} +\bigl(-\frac{39}{578} i b_{0}^{7} + \left(-\frac{91}{578} i + \frac{91}{578}\right) b_{0}^{6} - \frac{649}{578} b_{0}^{5} + \left(-\frac{1499}{578} i - \frac{1499}{578}\right) b_{0}^{4} + \frac{2233}{578} i b_{0}^{3} + \left(\frac{2666}{289} i - \frac{2666}{289}\right) b_{0}^{2} + \frac{219}{578} b_{0} + \frac{861}{289} i + \frac{861}{289}\bigr)q^{3}+\left(b_{0}^{2} - 4 i\right)q^{4} $\newline$+\bigl(-\frac{15}{578} b_{0}^{7} + \left(-\frac{35}{578} i - \frac{35}{578}\right) b_{0}^{6} +\frac{383}{578} i b_{0}^{5} + \left(\frac{621}{578} i - \frac{621}{578}\right) b_{0}^{4} +\frac{2993}{578} b_{0}^{3} + \left(\frac{1181}{289} i + \frac{1181}{289}\right) b_{0}^{2} - \frac{6131}{578} i b_{0} - \frac{220}{289} i + \frac{220}{289}\bigr)q^{5} + O(q^{6}),$ for $\Gamma_1(53)$ with nebentype mod 53 mapping $2 \mapsto i$, with coefficients in $\mathbb{Q}(i)[b_{0}]/(b_{0}^{8} + \left(3 i + 3\right) b_{0}^{7} - 16 i b_{0}^{6} + \left(-52 i + 52\right) b_{0}^{5} - 48 b_{0}^{4} + \left(-207 i - 207\right) b_{0}^{3} - 26 i b_{0}^{2} + \left(122 i - 122\right) b_{0} - 7)$. \\
\hline
\pagebreak[2]\hline
\multicolumn{3}{|l|}{\textbf{Level 59}.  Dimension 36.} \\
\hline
III & 4 & $q + b_{1}q^{2} + \left(\frac{1}{4} b_{1}^{4} + \frac{9}{2} b_{1}^{2} + \frac{65}{4}\right)q^{3} + \left(b_{1}^{2} + 4\right)q^{4} + \left(\frac{1}{4} b_{1}^{4} + \frac{11}{2} b_{1}^{2} + \frac{93}{4}\right)q^{5} + O(q^{6}),$ for $\Gamma_1(59)$ with nebentype mod 59 mapping $2 \mapsto -1$, with coefficients in $\mathbb{Q}[b_{1}]/(b_{1}^{6} + 27 b_{1}^{4} + 215 b_{1}^{2} + 509)$. \\
\hline
\end{longtable}
\end{center}


\bibliographystyle{amsalpha_no_mr}
\bibliography{AGM-VI}

\end{document}